\documentstyle[12pt]{article}      
\input amssym.def
\input amssym.tex                                           
\newtheorem{theorem}{Theorem}

\title{Two Remarks on K\"{a}hler Homogeneous Manifolds}

\begin{document}  

%\maketitle

\begin{center}{\bf\Large Two Remarks on \\ 
K\"{a}hler Homogeneous Manifolds}
\end{center}  

\begin{center}
{Bruce Gilligan\footnote{Dept. of Mathematics and Statistics, 
University of Regina, Regina, Canada S4S 0A2 \\ 
Email: gilligan@math.uregina.ca}         
and Karl Oeljeklaus\footnote{Centre de Math\'{e}matiques et d'Informatique, CNRS-UMR 6632 (LATP), 39, rue Joliot-Curie, 
Universit\'{e} de Provence, 13453 Marseille Cedex 13 France \\ 
Email: Karl.Oeljeklaus@cmi.univ-mrs.fr}}  
\end{center}

%\hline 

\begin{abstract} 
We prove that every K\"{a}hler solvmanifold has a finite covering  
whose holomorphic reduction is a principal bundle.  
An example is given that illustrates the necessity, in general,  
of passing to a proper covering.  
We also answer a stronger version of a question posed by Akhiezer 
for homogeneous spaces of nonsolvable algebraic groups 
in the case where the isotropy has the property that its intersection with the radical is 
Zariski dense in the radical.  
\end{abstract}  

\begin{quote} 
\begin{center}{\bf R\'{e}sum\'{e}}  
\end{center} 
Nous montrons que chaque vari\'{e}t\'{e} homog\`{e}ne r\'{e}soluble 
admet une rev\^{e}tement finie dont la r\'{e}duction holomorphe est un fibr\'{e} principal.  
On donne un exemple qui montre qu'il est n\'{e}cessaire en g\'{e}n\'{e}ral
de passer \`{a} une rev\^{e}tement finie.  
Nous donnons de plus  une r\'{e}ponse \`{a} une version plus forte 
d'une question d'Akhiezer 
pour un espace homog\`{e}ne d'un groupe alg\'{e}brique non-r\'{e}soluble 
dans le cas o\`u l'isotropie est tel que l'intersection avec le radical 
est dense au sens de Zariski dans le radical.   
\end{quote}  

%\hline  

\section{Introduction}  
Let $G$ be a complex Lie group, $H \subset G$ a closed
complex subgroup and $X:=G/H$ the corresponding complex homogeneous manifold. 
One main question in the theory of complex homogeneous manifolds is to specify conditions 
on $G$ and $H$ under which the manifold $X$ admits a (not
necessarily $G$-invariant) K\"ahler metric.   
For $X$ compact the theorem
of Borel-Remmert (cf. \cite{BorRem}) gives a complete answer to the problem.
In the non-compact case answers are known under certain assumptions on $G$,
e.g. $G$ semisimple \cite{Bert} or solvable \cite{OeRi1}. 
In this paper we make two further contributions. 

We first consider $X=G/H \to G/I =Y$ the holomorphic reduction of $X$, i.e., 
the unique $G$-equivariant fibration such that all holomorphic functions
on $X$ come from the holomorphically separable base $Y$, see \cite{GiHu}.
It has the universality property that every holomorphic map of $X$ into a 
holomorphically separable complex manifold 
factors through this fibration. 
Using known results, see \cite{OeRi1}, we prove that if $G$ is solvable
and $X$ is K\"ahler then the holomorphic reduction of $X$ is, up to
a finite covering, a {\it principal bundle} with a Cousin group as fiber.
A {\it Cousin group} is a complex abelian Lie group admitting no non constant
holomorphic functions.  

Finally, we consider the following stronger version of a question posed by Akhiezer \cite{Akh1}:  
{\it Suppose $G$ is a 
connected complex Lie group, $H$ a closed complex subgroup of $G$ and $X=G/H$
with ${\cal O}(X) \cong {\Bbb C}$ and there is no $G$-equivariant fibration 
of $X$ over a positive dimensional homogeneous rational manifold.
Is it true that $X$ K\"ahler implies that $X$ is a Cousin group?
 }  \\
We prove that if $G$ is non-solvable algebraic and if the intersection of
the radical $R \subset G$ with $H$ is Zariski dense in $R$ then $X$ is not K\"ahler.

\section{K\"ahler solvmanifolds}  

Suppose $X=G/H$ is a K\"{a}hler solvmanifold, i.e., $G$ is a complex solvable Lie group.  
We consider its holomorphic reduction 
$$  
  G/H \; \longrightarrow \; G/I  %\leqno{(*)}
$$
and prove that it is always a principal bundle up to a finite covering.  
It is already known that its base $G/I$   
is Stein \cite{HuOe} and that its fiber $C :=I/H$ is connected and is, 
as a complex manifold, biholomorphic to the Cousin group $I^{\circ}/H \cap I^{\circ}$  
\cite{OeRi1}. 

\begin{theorem}  
Suppose $G$ is a connected, solvable, 
complex Lie group and 
$H$ is a closed complex subgroup of $G$ with $X := G/H$ a
K\"{a}hler manifold.  
Let 
$$  
    G/H \; \longrightarrow \; G/I  %\leqno{(*)}
$$
be the holomorphic reduction.  
Then there is a subgroup of finite index $\hat I \subset  I$ such that the bundle 
$$ 
   \hat X:= G/{\hat I \cap H} \; \longrightarrow \; G/{\hat I}  %\leqno{(*)}
$$ is a holomorphic $I^{\circ}/H \cap I^{\circ}$-principal bundle and represents the
holomorphic reduction of $\hat X$.  
\end{theorem}  

\noindent{\bf Proof:}\  
Without loss of generality we may assume that $G$ is simply connected.
Then $G$ admits a Hochschild-Mostow 
hull, see \cite{HochMost},  i.e., there exists a 
solvable linear algebraic group 
\[ 
     G_{a} \; = \;  (\Bbb{C}^{*})^{k} \ltimes G  
\] 
and $G_{a}$ contains $G$ as a Zariski dense, topologically closed, 
normal complex subgroup.  

\noindent{\bf Part I:}\  We first assume that the isotropy is discrete.
In order to underline this fact, we will write the 
isotropy subgroup as $\Gamma$ in this part of the proof.  
Let 
$$ 
   X \; = \; G/\Gamma \; \longrightarrow \; G/J\cdot\Gamma 
$$ 
be the holomorphic reduction of $X$ with $J$ connected.  
Because $X$ is K\"{a}hler, the fiber $J\cdot\Gamma/\Gamma$ 
of its holomorphic reduction is (as a complex manifold) isomorphic to a Cousin group 
\cite{OeRi1}.  
We remark now that $G_{a}/\Gamma$ is biholomorphic to 
$(\Bbb{C}^{*})^{k}\times (G/\Gamma)$.  
As a consequence, $G_{a}/\Gamma$ is K\"{a}hler and 
its holomorphic reduction is given by 
$$    
     G_{a}/\Gamma \; \longrightarrow \; 
     G_{a}/J\cdot\Gamma \; \cong \; (\Bbb{C}^{*})^{k} \times 
     (G/J \cdot\Gamma ) , 
$$       
because the base of this fibration is holomorphically separable and its fiber is 
$J\cdot\Gamma/\Gamma$ and is, as a complex manifold, biholomorphic to a Cousin group.  
Let $G_{a}(\Gamma)$ denote the algebraic closure of $\Gamma$ in
$G_{a}$.  
The manifold $G_{a}/G_{a}(\Gamma)$ is Stein,
as a quotient of algebraic solvable Lie groups. 
In view of the universality property of the holomorphic reduction one has the following commutative diagram  
$$ 
         \begin{array}{ccc}  G_{a}/\Gamma & \longrightarrow & G_{a} / J\cdot \Gamma \\  
           &  \searrow  & \downarrow \\  
             & & G_{a} / G_{a}(\Gamma)    \end{array}   
$$
In particular, $J\cdot\Gamma \subset G_{a}(\Gamma)$.

Now consider the fibration 
$$ 
    G_{a}(\Gamma)/\Gamma \; \longrightarrow \;
     G_{a}(\Gamma)/J\cdot\Gamma  . 
$$  
This fibration is the holomorphic reduction of $G_{a}(\Gamma)/\Gamma$, 
because its fiber is, as a complex manifold a Cousin group, and its base, as an orbit in 
$G_{a}/J\cdot\Gamma$, is holomorphically separable and hence Stein.  \\
Let  $G_{a}(\Gamma)^{\circ}$ be the connected component of the identity
of $G_{a}(\Gamma)$
and define $\hat \Gamma:= \Gamma\cap G_{a}(\Gamma)^{\circ} $.
Then $\hat \Gamma$ is a 
Zariski dense subgroup of $G_{a}(\Gamma)^{\circ}$. 
Because the quotient $G_{a}(\Gamma)^{\circ}/\hat \Gamma $ is K\"{a}hler 
and $\hat \Gamma$ is Zariski dense in $G_{a}(\Gamma)^{\circ}$, it follows that
$G_{a}(\Gamma)^{\circ}$ is a unipotent,
see Lemma 3, p. 412 in \cite{OeRi1}.  

Therefore 
$$ 
    G_{a}(\Gamma)^{\circ}/\hat \Gamma \; \longrightarrow \;
G_{a}(\Gamma)^{\circ}/J\cdot\hat \Gamma  \ \ \leqno{(\dagger)}
$$  
is the holomorphic reduction of the homogeneous K\"ahler nilmanifold
$G_{a}(\Gamma)^{\circ}/\hat \Gamma$ the fiber
of which is still biholomorphic to  $J / \Gamma \cap J$.
  
Let us introduce the following notation:  
Let ${\frak g}_{a}$ denote the Lie algebra of the connected unipotent group 
$G_{a}(\Gamma)^{\circ}$, and $Z$ denote its center, 
with its Lie algebra being ${\frak z}$. 
Since the exponential map $\exp: \frak g \to G_{a}(\Gamma)^{\circ}$
is bijective we may consider the real Lie algebra ${\frak g}_{\Bbb R}$ spanned by 
$\exp^{-1}(\hat \Gamma)$ and its maximal complex ideal  
${\frak m} :=  {\frak g}_{\Bbb R}\cap i {\frak g}_{\Bbb R}$.   
 
The construction of the holomorphic reduction 
of a complex nilmanifold given in \cite{GiHu} shows that 
the group $J$ is the smallest connected, closed, complex subgroup 
of $G_{a}(\Gamma)^{\circ}$ that contains the connected group 
corresponding to the ideal ${\frak m}$ 
and such that $J\cdot\hat{\Gamma}$ is closed in $G_{a}(\Gamma)^{\circ}$.  
We note that $Z$ is a group such that $Z \cdot\hat{\Gamma}$ is closed, as 
follows from a result of Barth-Otte; see Theorem 4 in \cite{GiHu}.    
By using the K\"ahler assumption and detailed Lie algebra calculations  
it is shown in the proof of Theorem 2',  p. 409 in \cite{OeRi1}  that 
$$  
         {\frak m} \; \subset \; {\frak z} .  
$$ 
As a consequence, $J\subset Z$.  
In particular, $J$ centralizes $\hat \Gamma$ in 
the nilpotent group $G_{a}(\Gamma)^{\circ}$.  
Hence $ \hat \Gamma  \lhd J\cdot  \hat \Gamma$ and thus 
$J\cdot  \hat \Gamma/ \hat \Gamma$ is a group.  We note also that $\hat I:=J \cdot  \hat \Gamma$ is a 
subgroup of finite index of $J \cdot \Gamma$ and $\hat{I} \cap \Gamma =  \hat{\Gamma}$.
Thus the holomorphic reduction $(\dagger)$       
is a $J\cdot  \hat \Gamma/  \hat \Gamma=\hat I/\hat \Gamma$-principal bundle and 
the theorem is proved in this case.
 
\vskip 2ex\noindent{\bf Part II:}\  
We now assume that the isotropy 
subgroup $H$ is not discrete in $G$.  
Let 
$$  
    G/H \; \longrightarrow \; G/I 
$$
be the holomorphic reduction of $G/H$.  
We would like to show that $H$ is a normal subgroup of $I$.  
Here again $G_{a}/H$ is biholomorphic to 
$({\Bbb C}^{*})^{k} \times G/H$ and is therefore K\"{a}hler.      
Also $G_{a}/H \to G_{a}/I$ is clearly the 
holomorphic reduction of $G_{a}/H$.  
Set $N := (N_{G_{a}}(H^{\circ}))^{H}$, i.e., those connected
components of the normalizer 
of $H^{\circ}$ in $G_{a}$ that meet $H$.  
Note that $N/H$ is connected and  
$N$ is an algebraic subgroup of $G_{a}$.   
Therefore, $N^{\circ}$ has finite index in $N$.  
We remark that $G_{a}/N$ is holomorphically separable by Lie's Flag Theorem.    
Again by the universality property of the holomorphic reduction one has the following diagram  
$$ 
         \begin{array}{ccc}  G_{a}/H & \longrightarrow & G_{a} / I \\  
           &  \searrow  & \downarrow \\  
             & & G_{a} / N  \end{array}   
$$  
In particular, $I \subset N$.  
Since $N/I$ is an orbit in $G_{a}/I$ and the latter space is 
Stein, it follows that $N/I$ is holomorphically separable.   
But we are in a K\"{a}hler setting and thus the fiber $I/H$ 
of the holomorphic reduction of $G/H$ is 
biholomorphic to a (connected) Cousin group.  
As a consequence, $N/H \to N/I$ is the holomorphic
reduction of $N/H$.  Now define $\hat I := I \cap N^{\circ}$.
Then $\hat I \cap H = N^{\circ} \cap H$. Moreover,
$\hat I \subset G$ and $\hat I \cap H \subset G$.
With this notation we get as in the first case the fibration 
$$N/H \simeq
N^{\circ}/(  N^{\circ} \cap H) 
\to N^{\circ}/\hat I \simeq N/I$$ 
which is the holomorphic reduction of $
N^{\circ}/( N^{\circ} \cap H) $ with fiber $\hat I/(  N^{\circ} \cap H)$.

Now we have the following 
$$
        ( N^{\circ} / H^{\circ} ) / (  N^{\circ} \cap H /H^{\circ}) \; \simeq \;   N^{\circ} / 
          N^{\circ}\cap H  \; \longrightarrow \;   N^{\circ} / \hat I  \; \simeq \;            
      ( N^{\circ} / H^{\circ} ) / (\hat I /H^{\circ})       
$$ 
and Part I above applies.   
It follows that $ N^{\circ} \cap H/H^{\circ}$ is a normal
subgroup of $( N^{\circ} \cap H/H^{\circ})\cdot(I^{\circ}/H^{\circ})$  
and that $ N^{\circ}\cap H$ is a normal subgroup of $\hat I$ which
in turn is a subgroup of finite index of $I$.  
Thus the holomorphic reduction of the K\"{a}hler solvmanifold $G/H$ is, 
up to a finite covering, a Cousin group principal bundle.  
\hfill $\Box$

\medskip

The following examples show that, in general, passing to a subgroup of finite index is
unavoidable in order to get a principal bundle and that the K\"ahler condition is necessary.

\medskip

\noindent   
{\bf Example 1:}\
Let
$$  
   G \;  := \;  \left\{ \left(
                \begin{array}{ccc}
                a & b     \\
                0 & 1  \end{array}
                \right)  \mid \ a \in {\Bbb C}^*,\ b \in {\Bbb C} \ \right\}  
$$
and
$$
        \Gamma \; := \; \left\{ \left(
              \begin{array}{ccc}
              \sqrt{-1}^{\, n} & k + \sqrt{-1} l     \\
              0 &    1  
              \end{array}
               \right)  \mid \  n,k,l \in \Bbb Z \ \right\} .  
$$     
Then $X:= G/\Gamma$ is K\"ahler, $\Gamma$ is {\it not} nilpotent,
and the holomorphic reduction is {\it not} a principal bundle.
But of course all this is true, if one replaces $\Gamma$
by the obvious subgroup of index $4$.

\medskip\noindent 
\noindent   
{\bf Example 2:}\
We also note that the K\"{a}hler condition is necessary.     
We give the following example due to J.J. Loeb \cite{Loeb}.
Let $ K =\Bbb Z, \Bbb C$ and 
$G_K = K \ltimes K^2$ be defined by: 
$$ (z,b)\circ (z',b'):= (z+z', e^{Az'} b +b'),$$
where $z,z' \in K, \ \ b,b' \in K^2$ and $A$ is the real logarithm
of the matrix 
$\tiny
(\! \!
\begin{array}{cc}
2 \! &\! 1     \\
 1 \! &\!  1    
\end{array}
\! \!)
$. 
The complex solvmanifold $X= G_{\Bbb C}/G_{\Bbb Z}$ is not K\"ahler and has as holomorphic reduction  
$$  
       G_{\Bbb C}/G_{\Bbb Z}  \; \stackrel{{\Bbb C}^{*} \times {\Bbb C}^{*}}{\longrightarrow} \; 
       G_{\Bbb C}/ G_{\Bbb C}'\cdot G_{\Bbb Z} \; = \; {\Bbb C}^{*}  
$$  
No finite covering of $X=G_{\Bbb C}/G_{\Bbb Z}$ has a holomorphic reduction 
that is a principal bundle since in this case $X$ would be Stein and thus K\"{a}hler.

\section{A non-K\"ahler criterion}
 
For a complex Lie group $G$ we let $G=S \cdot R$ be a Levi-Malcev
decomposition, where $S$ is semisimple and $R$ is the radical of $G$.

\begin{theorem}  
Suppose $G=S\cdot R$ is a Levi-Malcev decomposition of an algebraic group with both 
factors having positive dimension.  
Let $H$ be a closed complex subgroup of $G$ and $X:=G/H$ such that  
${\mathcal O}(X) \cong {\Bbb C}$.   
Assume $H$ is not contained in any proper parabolic subgroup of $G$,
i.e.  there is no $G$-equivariant fibration 
of $X$ over a positive dimensional homogeneous rational manifold.
Suppose $R\cap H$ is Zariski dense in $R$.   
Then $G/H$ is not K\"{a}hler.  
\end{theorem}

\noindent{\bf Proof:}\  
Without loss of generality we may assume that the $G$-action on $X$
is almost effective. 
Let $N:=N_G(H^{\circ})$ be the normalizer of the identity component
of $H$ in $G$. 
Then $G/N$ is an orbit in some projective space and our assumptions force
$N=G$, i.e. $H$ is discrete. 
We shall use $\Gamma$ instead of $H$ for this reason.
We assume that $X$ is K\"{a}hler and derive a contradiction.  
It is known that the orbits of $R$ are closed in our situation, see \cite{OeRi2}, 
and so we may consider the fibration 
$$ 
      G/\Gamma \; \longrightarrow \;  G/R\cdot\Gamma \; = \; S/(S \cap R\cdot\Gamma)  .  
$$  
Note that $\Lambda := S \cap R\cdot\Gamma$ is a discrete subgroup of $S$.  
It follows from the assumptions on $X$ that $\Lambda$ is a Zariski dense subgroup of $S$,
since every algebraic subgroup of a semisimple complex Lie group is contained
in either a maximal reductive or a maximal parabolic subgroup of $G$.
Recall that a homogeneous manifold of a semisimple Lie group is Stein
if and only if the isotropy group is reductive, \cite{Mat}, \cite{On}.\\
As a consequence, we can find a semisimple element $\lambda\in\Lambda$
of infinite order.  
The Zariski closure of the infinite cyclic subgroup generated by $\lambda$ 
is an algebraic torus $A := \overline{\langle \lambda \rangle_{\Bbb Z}}
\cong ({\Bbb C}^{*})^{l}$ for some $l>0$.  
Without loss of generality, we may assume $A$ is connected, since one can 
replace $\lambda$ by a finite power, if necessary. 
Note that $A/\langle \lambda \rangle$ is a Cousin group and therefore
$A/A \cap \Lambda$ is also Cousin.   
Now let $S = S_{1} \cdots S_{n}$ be a decomposition of $S$ into 
a product of simple factors and let $p_{j}:S \to S_{j}$ denote the 
projection onto the $j$-th factor $S_{j}$.  
It was noted in the Lemma on p. 328 in \cite{Akh} that one may choose 
$\lambda$ so that 
$$p_{j}(A) \not= \{ e \} {\rm \  for \  all }\  j=1, \ldots , n. \leqno{(*)}  $$
 
Set $G^{*} := A \ltimes R$ which is an algebraic solvable subgroup of $G$ and
$\Gamma^{*} := G^{*}\cap \Gamma $.  
Let $N$ be the Zariski closure of $\Gamma^{*}$ in $G^{*}$.    
Then $N \subset G^*$ and $N^{\circ}$ is nilpotent by Lemma 4.3, p. 412, in \cite{OeRi1}.  

We claim $N=G^{*}$.  
Consider the fibration 
$\pi:G^{*}/\Gamma^{*} \to A/A \cap \Lambda$ 
and its restriction to the $N^{\circ}$-orbit $N^{\circ}/\Gamma^{*}\cap N^{\circ}$.  
The fiber of the restriction $\pi |_{N^{\circ}/\Gamma^{*}\cap N^{\circ}}$ is 
equal to $R/(R\cap\Gamma)$ 
by our assumption that $\Gamma \cap R$ is Zariski dense in $R$.  
This implies $R \subset N^{\circ}$.  
Now the projection $p:G \to S \cong G/R$ is an 
algebraic homomorphism and thus the image of an algebraic 
subgroup of $G$ containing $R$ is an algebraic subgroup of $S$.    
Hence $p(N^{\circ})$ is an algebraic subgroup of $S$ that contains 
$\langle\lambda\rangle_{\Bbb Z}$.  
As a consequence, $p(N^{\circ})=A$.  
These two observations imply $N^{\circ}=G^{*}$. 
Thus $G^{*}$ is nilpotent.  

Let $\phi$ denote the representation of $A$ on $R$ 
in the semidirect product $G^{*}= A\ltimes_{\phi}R$.  
Since $A \cong({\Bbb C}^{*})^{l}$, it follows that 
$\phi(A)$ is reductive.  
But $\phi(A)$ is also nilpotent and so must be trivial.    
Thus $A$ acts trivially on the radical $R$.  
Therefore $G^{*}= A \times R$ is a direct product.  
But this implies that the group $G$ is a group theoretic direct product 
$S\times R$ by the choice of $A$, see $(*)$, and the fact, just proved, that  
the subgroup $A$ is contained in the kernel of the representation 
that builds the semidirect product $S\ltimes R$.  
This yields the desired contradiction, 
since $S$ cannot have positive dimension in the case of a K\"ahler quotient
of a direct product $S \times R$, see \cite{OeRi2}.  
\hfill $\Box$

\vskip 2ex\noindent{\bf Acknowledgements:}  The first author gratefully acknowledges  
that he was partially supported by an NSERC Discovery Grant and was an 
Invit\'{e} du LATP \`{a} l'Universit\'{e} de Provence.

\end{document}